\documentclass[11pt]{article}
\usepackage[scale=0.9]{geometry}
\usepackage{graphics}
\usepackage{amsmath, amsthm,amssymb}
\usepackage{delarray}
\usepackage{leftidx}
\makeatletter
\newcommand*{\rom}[1]{\expandafter\@slowromancap\romannumeral #1@}
\makeatother

\def\1{\mathbf{1}}

\ifx\pdfoutput\undefined
\usepackage[pdftex]{graphicx}
\else
\usepackage[dvips]{graphicx}
\fi
\usepackage{amsmath,empheq}
\usepackage{xcolor}

\usepackage{listings}

\newtheorem{defn} {Definition}

\newtheorem{rem} {Remark}

\begin{document}

\sloppy

\date{20/11/2011}

\title{
\author{Controlled Continuous Time Random Walks\\
and Fractional Hamilton Jacobi Bellman Equations \footnotemark \\
\quad \nonumber \\
V. Kolokoltsov, M. Veretennikova \footnotemark}}

\sloppy

\date{March 26th, 2012}

\maketitle

\vspace{-0.2cm}

\footnotetext[1]{University of Warwick, Mathematics Institute, CV4 7AL}
\footnotetext[2]{Supported by EPSRC and MASDOC DTC, grant EP/HO23364/1}

\smallskip

\begin{abstract}
In this paper we study controlled continuous time random walks (CTRWs) and heuristically derive pay-off function dynamic programming (DP) equations which turn in the limit of standard scaling to fractional Hamilton Jacobi Bellman (fHJB) type equations. This paper aims to extend the results from \cite{kolokoltsov2009generalized} in a controlled setting.
\end{abstract}

\smallskip


\section{Introduction}

\subsection{Plan}

We introduce the controlled process and write the DP equation for the pay-off function. We scale the waiting times and the jump sizes by a parameter $\tau$ raised to specified powers, and obtain the limiting equation for the optimal pay-off function as $\tau \rightarrow 0$. Then we make slight modifications to the process, such as including waiting and jumping reward functions as well as an inner motion during waiting, and write corresponding DP equations for these pay-off functions, again looking at what happens as $\tau \rightarrow 0$. After this we study a generalised version of such a process and arrive at a general limiting equation for the pay-off function.

\subsection{Notation}

\begin{itemize}
\item $\leftidx{_{a+}}{D^{\beta}_{x}}{}$ left-sided Riemann-Liouville fractional derivative of order $\beta \in (0,1)$, with a lower integration limit $a$, which can be $-\infty$

\item $\leftidx{_{b-}}{D^{\beta}_{x}}{}$ right-sided Riemann-Liouville fractional derivative of order $\beta \in (0,1)$, with a lower integration limit $b$, which can be $\infty$

\item $\leftidx{_{a+}}{D^{* \beta}_x}{}$ left-sided Caputo fractional derivative of order $\beta \in (0,1)$, with a lower integration limit $a$, which can be $-\infty$

\item $\leftidx{_{b-}}{D^{* \beta}_{x}}{}$ right-sided Caputo fractional derivative of order $\beta \in (0,1)$, with a lower integration limit $b$, which can be $\infty$

\item $\frac{d^{\beta}}{d y^{\beta}}$ the generator form of the fractional derivative of order $\beta \in (0,1)$

\item $[x]$ is the integer part of $x \in \mathbb{R}$

\item For a complete metric space $S$, $C_{\infty}(S)$ consists of functions $f$ such that $\lim_{x \rightarrow 0}f(x)=0$, i.e. $\forall \epsilon > 0$, there exists a compact set $K$ such that $\sup_{x \notin K}|f(x)| < \epsilon$

\item $S(\mathbb{R}^{d})=\{f \in C^{\infty}(\mathbb{R}^{d}: \forall k,l \in \mathbb{N}, |x|^{k}\nabla^{l}f(x) \in C_{\infty}(\mathbb{R}^{d}))\}$ is the Schwartz space of rapidly decreasing functions

\item $C^{k}(\mathbb{R}^{d})$ is the Banach space of $k$ times continuously differentiable functions with bounded derivatives on $\mathbb{R}^{d}$, the norm is the sum of the sup-norms of the function and all of its derivatives of orders up to $k$ inclusively

\item $\mathbf{1}_{\{z \in A\}}$ stands for the indicator function of a set $A$, which is equal to 1 when its argument $z$ is in the set $A$, otherwise it is $0$
\end{itemize}

\section{Fractional calculus}

There are numerous articles and books with an introduction to fractional calculus and its uses, see e.g. \cite{meerschaert2012stochastic}, \cite{podlubny1999fractional} and \cite{kilbas2006theory}. Let's recall the main definitions and links for convenience of the reader. The Caputo derivative is one of the several versions of the fractional derivative. It is closely interlinked with the generator form of the fractional derivative and the Riemann Liouville derivative in the way described below.

The left-sided Caputo derivative for $\beta \in (0,1)$ is defined for $f$ smooth on $[a,\infty)$ by
\begin{equation}
\leftidx{_{a+}}{D^{* \beta}}{_x}f(x)=\frac{1}{\Gamma(1-\beta)}\int_{a}^{x}\frac{d f(y)}{dy}(x-y)^{-\beta}dy.
\end{equation}

The right-sided Caputo derivative for $\beta \in (0,1)$ is defined for $f$ smooth on $(-\infty,b]$ by
\begin{equation}
\leftidx{_{b-}}{D^{* \beta}}{_x}f(x)=-\frac{1}{\Gamma(1-\beta)}\int_{x}^{b}\frac{d f(y)}{dy}(y-x)^{-\beta}dy.
\end{equation}

The left-sided Riemann-Liouville fractional derivative for $\beta \in (0,1)$ is defined for $f$ smooth on $[a,\infty)$ as follows:
\begin{equation}
\leftidx{_{a+}}{D^{\beta}}{_x}f(x)=\frac{1}{\Gamma(1-\beta)}\frac{d}{dx}\int_{a}^{x}f(y)(x-y)^{-\beta}dy.
\end{equation}

The right-sided Riemann-Liouville fractional derivative for $\beta \in (0,1)$ is defined for $f$ smooth on $(-\infty,b]$ as follows:
\begin{equation}
\leftidx{_{b-}}{D^{\beta}}{_x}f(x)=\frac{1}{\Gamma(1-\beta)}\frac{d}{dx}\int_{x}^{b}f(y)(y-x)^{-\beta}dy.
\end{equation}

The left-sided fractional derivatives operate on $a<y<x$, whilst the right-sided fractional derivatives operate on $x<y<b$.

For every function $f$ which is everywhere continuously differentiable \cite{kilbas2006theory},

\begin{equation}\label{RLCapa}
\leftidx{_{a+}}{D^{\beta}}{_x}f(x)=\leftidx{_{a+}}{D^{* \beta}}{_x}f(x)+ \frac{(x-a)^{-\beta}}{\Gamma(1-\beta)}f(a).
\end{equation}

To prove this do a change of variables $x-y=z$ and differentiate the integral, then let $z=x-y$ to get the result:

\begin{eqnarray}
\leftidx{_{a+}}{D^{\beta}}{_x}f(x)=\frac{1}{\Gamma(1-\beta)}\frac{d}{dx}\int_{a}^{x}f(y)(x-y)^{-\beta}dy\nonumber \\
=-\frac{1}{\Gamma(1-\beta)}\frac{d}{dx}\int_{x-a}^{0}f(x-z)z^{-\beta}dz \nonumber \\
=\frac{1}{\Gamma(1-\beta)}\frac{d}{dx}\int_{0}^{x-a}f(x-z)z^{-\beta}dz \nonumber \\
=\frac{1}{\Gamma(1-\beta)}\left[(x-a)^{-\beta}f(x-x+a) + \int_{0}^{x-a}\frac{\partial f(x-z)}{\partial x}z^{-\beta}dz  \right] \nonumber \\
=\frac{1}{\Gamma(1-\beta)}\left[(x-a)^{-\beta}f(a) + \int_{a}^{x}\frac{df(y)}{dy}(x-y)^{-\beta}dy  \right], \nonumber \\
\end{eqnarray}
Q.E.D. Similarly one can analogously prove that

\begin{equation}\label{RLCapb}
\leftidx{_{b-}}{D^{\beta}}{_x}f(x)=\leftidx{_{b-}}{D^{* \beta}}{_x}f(x)+\frac{(b-x)^{-\beta}}{\Gamma(1-\beta)}f(b).
\end{equation}

The equations (\ref{RLCapa}) and (\ref{RLCapb}) motivate the definition of the {it regularized} Caputo derivative. When $f$ is not differentiable the {\it regularized} Caputo derivative is defined as
\begin{equation}
\leftidx{_{a+}}{\partial^{* \beta}}{_x}f(x)=\frac{1}{\Gamma(1-\beta)}\frac{d}{dx}\left[\int_{a}^{x}f(x-y)y^{-\beta}dy\right]-f(a)\frac{(x-a)^{-\beta}}{\Gamma{(1-\beta)}}.
\end{equation}

The generator form of the  {\it{fractional differential operator}} for every $f$ at least from the Schwartz space $S(\mathbb{R})$ is defined as follows:
\begin{equation}\label{gen}
\frac{d^{\beta}f(x)}{d x^{\beta}}=\frac{1}{\Gamma(-\beta)}\int_{0}^{\infty}(f(x-y)-f(x))\frac{dy}{y^{1+\beta}},\quad \quad \mbox{$\beta \in (0,1)$}.
\end{equation}
This also applies to all $f$ which are bounded and $\gamma$-Holder continuous, for $\gamma < \beta$ \cite{kolokoltsov2011markov}. We can also let $f$ be such that $f(x)=0$ for $x < a$ and $C^{1}[a,\infty)$. The reason for which this fractional derivative is referred to as the generator form one, is that for a $\beta$-stable Levy motion the generator is exactly,

\begin{equation}\label{stablegene}
A_{\beta}f(x)=-\frac{1}{\Gamma(-\beta)}\int_{0}^{\infty}(f(x+y)-f(x))\frac{dy}{y^{1+\beta}}
\end{equation}
which knots fractional calculus with probability. 

The dual operator to $A_{\beta}$ is defined as follows:

\begin{equation}\label{stablegene}
A^{*}_{\beta}f(x)=-\frac{1}{\Gamma(-\beta)}\int_{0}^{\infty}(f(x-y)-f(x))\frac{dy}{y^{1+\beta}}.
\end{equation}

In particular

\begin{equation}
-\frac{1}{\Gamma(1-\beta)}t^{-\beta}=-A^{*}_{\beta}\mathbf{1}_{\{t>0\}}, 
\end{equation}
where the dual to A operator 
\begin{equation}
A^{*}_{\beta}=\leftidx{_{a+}}{D^{* \beta}_{x}}{}=\leftidx{_{a+}}{D^{\beta}_{x}}{}
\end{equation}
with $a=-\infty$.

For every $a <0$, if $f(x)=0$ for $x<0$, then the generator form coincides with the Riemann-Liouville fractional derivative and with the Caputo fractional derivative. In case $a=0$, $f(x) = 0$ for $x <0$, the generator form is equivalent to Caputo form, but not to the Riemann-Liouville form \cite{meerschaert2012stochastic}. To prove the generator form is equivalent to the Caputo form for $a \le 0$, and $f(x)=0$ if $x<0$, integration by parts yields \cite{meerschaert2012stochastic}

\begin{eqnarray}\label{GC}
\int_{0}^{\infty}[f(x-y)-f(x)]\frac{1}{\Gamma(-\beta)}y^{-\beta-1}dy\nonumber \\
=\frac{y^{-\beta}(f(x-y)-f(x))}{-\beta\Gamma(-\beta)}|_{0}^{\infty}\nonumber \\
+\frac{1}{\beta \Gamma(-\beta)}\int_{0}^{\infty}\frac{d}{dy}[f(x-y)-f(x)]y^{-\beta}dy.\\
\end{eqnarray}

For $f$ bounded and continuously differentiable, for $y \rightarrow 0$

\begin{equation}
y^{-\beta}[f(x-y)-f(x)] = O(y^{1-\beta}),
\end{equation}

and for $y \rightarrow \infty$

\begin{equation}
y^{-\beta}[f(x-y)-f(x)] = O(y^{-\beta}).
\end{equation}

Hence the boundary term in (\ref{GC}) is $0$ and
\begin{eqnarray}
\frac{1}{\beta \Gamma(-\beta)}\int_{0}^{\infty}\frac{d}{dy}[f(x-y)-f(x)]y^{-\beta}dy\nonumber \\
=\frac{1}{\Gamma{(1-\beta)}}\int_{0}^{\infty}f'(x-y)y^{-\beta}dy \nonumber \\
=\frac{1}{\Gamma{(1-\beta)}}\int_{-\infty}^{x}\frac{df(y)}{dy}(x-y)^{-\beta}dy\nonumber \\
=\leftidx{_{(-\infty)+}}{D^{\beta}}{_x}f(x)\nonumber \\
=\leftidx{_{a+}}{D^{\beta}}{_x}f(x), \quad \mbox{$\forall a \le 0$}.\\
\end{eqnarray}




\section{A continuous time random walk (CTRW)}

For convenience we recall here the main definition of a CTRW.

\begin{defn}
For $i \in \mathbb{N}$, waiting times $\gamma_{i} \in \mathbb{R}_{+}$, and jump sizes $\xi_{i} \in \mathbb{R}^{d}$ we let $(\gamma_{1},\xi_{1})$, $(\gamma_{2},\xi_{2})$, $\ldots$ be a sequence of i.i.d. pairs of random variables (r.v.'s) such that the distribution of $\gamma_{i}$ is given by the probability measure $\nu(dr)$, and the distribution of jump sizes is given by the probability measure $\mu(d\xi)$. Let $X(n)= \sum_{i=1}^{n}\gamma_{i}$. Define the inverse process $Z_{X}(t)=\inf_{n}\{n: X(n) \ge t\}$, which is continuous in time. A {\it{continuous time random walk}} (CTRW) with jumps $\xi_{i} \in \mathbb{R}^{d}$ and waiting times $\gamma_{i} \in \mathbb{R}_{+}$ is defined as the process in continuous time given by 
\begin{equation}
Y(t)=Y_{Z_{X}(t)}=\sum_{i=1}^{Z_{X}(t)}\xi_{i}, 
\end{equation}
arising from the measures $\nu$ and $\mu$.
\end{defn}

One interprets this as follows: the particle waits for a time $\gamma_{i} \in \mathbb{R}_{+}$ and makes a jump $\xi_{i} \in \mathbb{R}^{d}$, then it waits for a time $\gamma_{i+1} \in \mathbb{R}_{+}$ and makes a jump $\xi_{i+1} \in \mathbb{R}^{d}$, e.t.c. for $i \in \mathbb{N}$. For now we assume that the waiting times and the jumps are not time or position-dependent random variables. Later in the paper, we do consider this complication too, see Chapter 8.

In a CTRW each jump is characterised by a waiting time and a random jump length that are i.i.d. random variables with probability densities $\nu$ and $\mu$. As it is usually done in the theory of CTRWs, e.g. \cite{meerschaert2012stochastic} and \cite{gorenflo2010mittag}, let us assume that for $i \in \mathbb{N}$, $\gamma_{i}$ and $\xi_{i}$ belong to domains of attraction of stable laws: for $n \rightarrow \infty$,
\begin{equation}\label{nustable}
\int_{|r| > n}\nu(dr) \sim \frac{1}{\Gamma(\beta-1)n^{\beta}},
\end{equation}
$\beta \in (0,1)$,
and concerning the jumps we assume for $n \rightarrow \infty$
\begin{equation}
\int_{|y|>n}\int_{\Omega} \sigma(\xi)d_{S}(\xi)\mu(dy) \sim \frac{1}{n^{\alpha}} \int_{\Omega}\sigma(\xi)d_{S}\xi,
\end{equation}
where $\alpha \in (0,2]$, $\sigma$ is a spectral measure on the sphere $S^{d-1}$,
and for every Borel $\Omega \in S^{d-1}$ which can be any sector of the unit sphere $S^{d-1}$ . If $\alpha=2$ we instead assume that 
\begin{equation}
\int_{\mathbb{R}^{d}} y^{2}\mu(dy)< \infty.
\end{equation}

\begin{rem}
Note that if $\gamma_{i}, i \in \mathbb{N}$ were exponentially distributed, i.e. $\beta=1$, then this is a standard random walk and the DP equation for the optimal pay-off will be the Hamilton-Jacobi-Bellman (HJB) equation with a classical time derivative $d/dt$. 
\end{rem}

\begin{rem}
As it is well known, see \cite{kolokoltsov2010mathematical}, \cite{kolmogorov1949sums}, the distributions of normalized sums of i.i.d. copies of such r.v.'s converge to a stable law.
\end{rem}



\section{Fractional Hamilton Jacobi Bellman (HJB) for a controlled scaled CTRW}

If we pay attention to the times when the jumps take place, then the process $Y(t)$ has Markovian properties at the jump times, which are the renewal times for the process studied only at jump times. It is possible that as the particle waits or jumps, it has to pay a price for doing so, or it receives reward. Here we will assume that we can control the distribution of jumps so that we can maximise our reward. We highlight that the jump control may be of any nature, and we denote by $U$ the set of all possible controls at any particular jump, the set of controls being the same for every jump. By $\tilde{U}$ denote the set of all possible controls at all the jump times:
\begin{equation}
\tilde{U}=\{\tilde{u}=(u_{1}, u_{2}, \ldots )\}, u_{i} \in U, i \in \mathbb{N}.
\end{equation}
Control brings about dependence of the probability measure $\mu$ on $u \in U$. 
When control is introduced into the system, for $i \in \mathbb{N}, \xi_{i}$ become dependent on control, so we write $\xi_{i}(u_{i})$, and the latter are control dependent and not identically distributed random variables. Their distribution now depends on control, hence the notation $\mu_{u}(d\xi)$ appearing later. Now 
\begin{equation}
Y(t,\tilde{u})=\sum_{i=1}^{Z_{X}(t)}\xi_{i}(u_{i})
\end{equation}
is our {\it{controlled CTRW}}. We are not going to control the distribution of waiting times, hence $Z_{X}(t)$ is the same as before.

Assume that the jumps are costless and pay-off is only the terminal one given by the function $S_{0}(y)$, and that there are no rewards after the process has terminated. Given $S_{0}(y)$, the optimal pay-off function $S(t,y)$, is then defined as follows
\begin{equation}
S(t,y)=\sup_{\tilde{u} \in \tilde{U}}\mathbb{E}S_{0}(y+Y(t,\tilde{u})), \quad t \ge 0.
\end{equation}

This represents the optimal reward the particle receives when the time until the end of the process is $t$ and the particle is in position $y$. Unless specified otherwise, we will refer to the time until the termination of the process as $t$. 

For convenience let us denote the waiting time until the next jump by the random variable $\gamma$ with distribution $\nu(dr)$ just as the previous waiting times $\gamma_{i}$, $i \in \mathbb{N}$, the waiting time until a jump $\xi$ with distribution $\mu_{u}(d\xi)$ occurs after time $t$. Then
\begin{eqnarray}
P[\gamma > t]=\int_{t}^{\infty}\nu(dr),
\end{eqnarray}
is the probability that there is no jump within time $t$ to the end of the process.
If a jump $\xi$ takes place at a time $r \in [0,t]$, then the expectation of the pay-off $S$ at $(t-r,y+\xi)$ will be represented by the integral $\int_{\mathbb{R}^{d}}S(t-r,y+\xi)\mu_{u}(d\xi)$,where $0<r<t$, $r=t-\gamma$.

Consequently, $S(t,y)$ satisfies the following DP equation:
\begin{eqnarray}\label{notscaled}
S(t,y)=\sup_{u \in U}\left[S_{0}(y)P[\gamma > t] + \int_{0}^{t}\mathbb{E} S(t-r,y + \xi(u)| \gamma = r)\nu(dr)\right] \nonumber \\
=\sup_{u \in U}\left[S_{0}(y)\int_{t}^{\infty}\nu(dr) + \int_{0}^{t}\int_{\mathbb{R}^{d}}S(t-r,y+\xi)\mu_{u}(d \xi) \nu(d r) \right].
\end{eqnarray}

The first term represents the expectation of the pay-off if a jump has not occurred within time $t$ to the end of the process, and the second term represents the expectation of the pay-off if a jump has occurred within that time.

Now let us make the following scaling: $\gamma_{i} \mapsto \gamma_{i} \tau^{1/\beta}$, and $\xi_{i} \mapsto \xi_{i} \tau^{1/\alpha}$, $i \in \mathbb{N}$, where $\tau > 0$ and is small, also $\beta \in (0,1)$ and $\alpha \in (0,2].$ For this scaled process, we will denote the optimal pay-off at time $t$ and position $y$ by $S^{\tau}(t,y)$:
\begin{equation}
S^{\tau}(t,y)=\sup_{\tilde{u} \in \tilde{U}}\mathbb{E}S^{\tau}_{0}(y+Y(t,\tilde{u})).
\end{equation}
Our previous $S(t,y)$ is $S^{1}(t,y).$ We are interested in the limit $\tau \rightarrow 0$. The DP equation for $S^{\tau}(t,y)$ becomes

\begin{eqnarray}\label{scaled}
S^{\tau}(t,y)=\sup_{u\in U}\left[S^{\tau}_{0}(y)P[\gamma > t]+ \int_{0}^{t}\mathbb{E} S^{\tau}(t-r,y+\xi(u))| \gamma=r)\nu(dr/\tau^{1/\beta})\right] \nonumber \\
=\sup_{u \in U}\left[S^{\tau}_{0}(y)\int_{t}^{\infty}\nu(dr/\tau^{1/\beta})+ \int_{0}^{t}\int_{\mathbb{R}^{d}}S^{\tau}(t-r,y+\xi)\mu_{u}(d\xi/\tau^{1/\alpha})\nu(dr/\tau^{1/\beta}) \right].
\end{eqnarray}
We assume that as $\tau \rightarrow 0$, $S^{\tau}$ has a limit $\hat{S}(t,y)$ which is smooth. This is the usual assumption of a heuristic approach to the deduction of differential equations of dynamic programming. We wish to identify the limiting equation satisfied by $\hat{S}(t,y)$.

\smallskip
By adding and subtracting identical terms, we can bring equation (\ref{scaled}) to the following form:

{\footnotesize
\begin{eqnarray}
S^{\tau}(t,y)=\sup_{u \in U}\left(S^{\tau}_{0}(y)\int_{t}^{\infty}\nu(dr/\tau^{1/\beta})+ \int_{0}^{t} \int_{\mathbb{R}^{d}} S^{\tau}(t-r,y+\xi)\mu_{u}(d\xi/\tau^{1/\alpha}) \nu(dr/\tau^{1/\beta})\right)\nonumber \\
=S^{\tau}_{0}(y)\int_{t \tau^{-1/\beta}}^{\infty}\nu(dr) + 
\sup_{u \in U} \int_{0}^{t \tau^{-1/\beta}} \int_{\mathbb{R}^{d}} 
\left(S^{\tau}(t-r \tau^{1/\beta},y+\xi \tau^{1/\alpha})-
S^{\tau}(t-r\tau^{1/\beta},y)+ S^{\tau}(t-r\tau^{1/\beta},y)\right) \mu_{u}(d\xi)\nu(dr)\nonumber \\
=\sup_{u \in U} \left[S^{\tau}_{0}(y)\int_{t \tau^{-1/\beta}}^{\infty}\nu(dr)+ \tau \int_{0}^{t \tau^{-1/\beta}} \int_{\mathbb{R}^{d}} \frac{1}{\tau} \left(S^{\tau}(t-r \tau^{1/\beta},y+\xi \tau^{1/\alpha})
-S^{\tau}(t-r\tau^{1/\beta},y) \right) \mu_{u}(d\xi)\nu(dr)\right. \nonumber \\ 
\left.+\tau \int_{0}^{t \tau^{-1/\beta}} \int_{\mathbb{R}^{d}} \frac{1}{\tau} 
[S^{\tau}(t-r\tau^{1/\beta},y)
-  S^{\tau}(t,y)] \mu_{u}(d\xi)\nu(dr) +  \int_{0}^{t \tau^{-1/\beta}} \int_{\mathbb{R}^{d}} S^{\tau}(t,y) \mu_{u}(d\xi)\nu(dr) \right] \nonumber \\
+\sup_{u \in U} \left[S^{\tau}_{0}(y)\int_{t \tau^{-1/\beta}}^{\infty}\nu(dr)+ \tau \int_{0}^{t \tau^{-1/\beta}} \int_{\mathbb{R}^{d}}\left[\frac{S^{\tau}(t-r\tau^{1/\beta},y+\xi \tau^{1/\alpha})-S^{\tau}(t-r\tau^{1/\beta},y)}{\tau} \right] \mu_{u}(d\xi)\nu(dr) \right. \nonumber \\
\left. +\tau \int_{0}^{t \tau^{-1/\beta}} \int_{\mathbb{R}^{d}} \left[ \frac{S^{\tau}(t-r\tau^{1/\beta},y)-S^{\tau}(t,y)}{\tau} \right]\mu_{u}(d\xi)\nu(dr) + \int_{0}^{t \tau^{-1/\beta}} \int_{\mathbb{R}^{d}} S^{\tau}(t,y) \mu_{u}(d\xi)\nu(dr)\right].\nonumber \\
\end{eqnarray}
}
Note that since $S^{\tau}(t,y)$ is independent of $\xi$, integration with respect to $\xi$ gives us $1$ in the inner integral because $\mu_{u}$ is a probability measure, and also since there is no dependence on control in the integral, the term can be taken out of the supremum. 

We divide the whole equation by $\tau$ and proceed as follows with the last term of the RHS:
\begin{eqnarray}
\int_{0}^{t \tau^{-1/\beta}} \int_{\mathbb{R}^{d}} S^{\tau}(t,y)\mu_{u}(d\xi)\nu(dr)\nonumber \\
=\int_{0}^{\infty} \int_{\mathbb{R}^{d}} S^{\tau}(t,y)\mu_{u}(d\xi)\nu(dr) - \int_{t \tau^{-1/\beta}}^{\infty} \int_{\mathbb{R}^{d}} S^{\tau}(t,y)\mu_{u}(d\xi)\nu(dr)\nonumber  \\
=\int_{0}^{\infty}  S^{\tau}(t,y)\nu(dr) - \int_{t \tau^{-1/\beta}}^{\infty}  S^{\tau}(t,y)\nu(dr).\nonumber  \\
\end{eqnarray}
We shall see both of these integral terms cancelling out with other terms in the equation. The term
$S^{\tau}(t,y)$ on LHS cancels with $\sup_{u \in U}\int_{0}^{\infty}S^{\tau}(t,y)\nu(dr)$ on the RHS, because $\nu$ is a probability measure and the latter integral is independent of control $U$.

Now we work with the second to last RHS term:
\begin{eqnarray}
\int_{0}^{t \tau^{-1/\beta}} \int_{\mathbb{R}^{d}} \left[ \frac{S^{\tau}(t-r\tau^{1/\beta},y)-S^{\tau}(t,y)}{\tau} \right]\mu_{u}(d\xi)\nu(dr) \nonumber \\
=\int_{0}^{t \tau^{-1/\beta}}  \left[ \frac{S^{\tau}(t-r\tau^{1/\beta},y)-S^{\tau}(t,y)}{\tau} \right]\nu(dr) \nonumber \\
=\int_{0}^{\infty}  \left[ \frac{S^{\tau}(t-r\tau^{1/\beta},y)-S^{\tau}(t,y)}{\tau} \right]\nu(dr)
  \nonumber \\
-\int_{t \tau^{-1/\beta}}^{\infty}  \left[ \frac{S^{\tau}(t-r\tau^{1/\beta},y)-S^{\tau}(t,y)}{\tau} \right]\nu(dr),
\end{eqnarray}
and in the second integral the first term is $0$ due to the boundary condition we have chosen.
The second term in the second integral cancels out with  $-\frac{1}{\tau} \int_{t \tau^{-1/\beta}}^{\infty}  S^{\tau}(t,y)\nu(dr)$, whilst the first integral, as $\tau \rightarrow 0$ gives us a neat expression presented just below. As $\tau \rightarrow 0$ we have \footnote{This convergence follows from our assumption on $\nu$ in (\ref{nustable}), and the details can be found for example in \cite{kolokoltsov2011markov}}
\begin{eqnarray}
\int_{0}^{\infty} \left[ \frac{S^{\tau}(t-r\tau^{1/\beta},y)-S^{\tau}(t,y)}{\tau} \right]\nu(dr) \rightarrow A^{*}_{\beta}\hat{S}(t,y),
\end{eqnarray}
where $A^{*}_{\beta}$ is the dual to the generator $A_{\beta}$ for $\hat{X}_{t}$, the limiting process of $X^{\tau}$ as $\tau \rightarrow 0$, where $A_{\beta}$ is defined in (\ref{stablegene}) and represents the generator of the $\beta$-stable Levy subordinator.

Now, we recall from the definition of the stable-like probability measure $\nu$:
\begin{equation}
\int_{n}^{\infty}\nu(dr) \sim \frac{1}{\Gamma(1-\beta)n^{\beta}}
\end{equation}
as $n \rightarrow \infty$. In our case $n=t \tau^{-1/\beta} \rightarrow \infty$ as $\tau \rightarrow 0$.

By the means of the above condition and the essential assumption that

\begin{empheq}[box={\color{black}\fboxsep=10pt\fbox}]{align}
\color{red}
S^{\tau}(t,y) \rightarrow \hat{S}(t,y), \mbox{ as } \tau \rightarrow 0,
\end{empheq}

where $\hat{S}(t,y)$ is smooth and belongs to the domains of the limiting stable generators $L_{\alpha}^{u}$ and $A_{\beta}$, the first term of the RHS, as $\tau \rightarrow 0$, gives us 
\begin{equation}
\lim_{\tau \rightarrow 0}\frac{S^{\tau}_{0}(y)}{\tau}\int_{t\tau^{-1/\beta}}^{\infty}\nu(dr) \sim \frac{\hat{S}_{0}(y)t^{-\beta}}{\Gamma(1-\beta)}.
\end{equation}

As for the other term on the RHS, as $\tau \rightarrow 0$ we have \footnotemark[1]

\begin{eqnarray}
\int_{0}^{t \tau^{-1/\beta}} \int_{\mathbb{R}^{d}}\left[\frac{S^{\tau}(t-r\tau^{\beta},y+\xi \tau^{1/\alpha})-S^{\tau}(t-r\tau^{1/\beta},y)}{\tau} \right] \mu_{u}(d\xi)\nu(dr) \rightarrow L_{\alpha}S^{\tau}(t-r\tau^{1/\beta},y) \rightarrow L_{\alpha}\hat{S}(t,y),
\end{eqnarray}
where $L_{\alpha}$ is the generator for $\hat{Y}_{t}$, the limiting process of $Y^{\tau}$ as $\tau \rightarrow 0$. Note that the stable like generator $L_{\alpha}$ has the following form for $\alpha > 1$:
\begin{equation}
L_{\alpha}f(y)= \int_{\mathbb{R}^{d}}\frac{f(y+r)-f(y)-f'(y)r}{|r|^{d+\alpha}}\omega(u)dr;
\end{equation}
and the general stable symmetric law generator for $\alpha > 1$ has the form in polar coordinates:
\begin{equation}
L_{\alpha}f(y)=\int_{0}^{\infty}d|r| \int_{S^{d-1}}\frac{f(y+r)-f(y)-(r,\nabla f(y))}{|r|^{\alpha + 1}}\omega(dS),
\end{equation}
where $\omega$ is an arbitrary symmetric finite Borel measure on $S^{d-1}.$

Whilst for $\alpha \in (0,1)$

\begin{equation}
L_{\alpha}f(y)= \int_{\mathbb{R}^{d}}\frac{f(y+r)-f(y)}{|r|^{d+\alpha}}\omega(u)dr;
\end{equation}
and the general stable symmetric law generator for $\alpha > 1$ has the form in polar coordinates:
\begin{equation}
L_{\alpha}f(y)=\int_{0}^{\infty}d|r| \int_{S^{d-1}}\frac{f(y+r)-f(y)}{|r|^{\alpha + 1}}\omega(dS),
\end{equation}

The above limiting procedures and that $\tau$ now cancels out in front of the integral terms of the RHS, now yield us:

\begin{center}
\begin{empheq}[box={\color{black}\fboxsep=10pt\fbox}]{align}\label{basic1}
\color{red}
0=\frac{1}{\Gamma(1-\beta)} t^{-\beta}\hat{S}_{0}(y)+\sup_{u \in U}\left[L_{\alpha}^{u}\hat{S}(t,y)\right] + A^{*}_{\beta}\hat{S}(t,y).
\end{empheq}
\end{center}

Here $A_{\beta}$ is the generator form of the fractional derivative and $A^{*}_{\beta}$ is its dual operator. Due to the link between the Riemann-Liouville and the generator form of the  fractional derivative when $a=0$ we can write this as:
\begin{center}
\begin{empheq}[box={\color{black}\fboxsep=10pt\fbox}]{align}\label{basic1}
\color{red}
\leftidx{_{0+}}{D^{ \beta}}{_x}\hat{S}(t,y)=\sup_{u \in U}\left[L_{\alpha}^{u}\hat{S}(t,y)\right].
\end{empheq}
\end{center}
recall that the latter fractional derivative $\leftidx{_{0}}{D^{ \beta}}{_x}$ is the Riemann-Liouville type. In case of no control this equation turns to

\begin{empheq}[box={\color{black}\fboxsep=10pt\fbox}]{align}
\color{red}
A^{*}_{\beta}\hat{S}(t,y)=-L_{\alpha}\hat{S}(t,y)+\hat{S}_{0}(y)A^{*}_{\beta}\mathbf{1}_{\{t>0\}}
\end{empheq}

which appears in \cite{kolokoltsov2011markov}.

The equation (\ref{basic1}) is fractional in time and in space, and is of Hamilton Jacobi Bellman (HJB) type. For controlled Markov diffusion processes instead of the fractional equation we have the second order HJB equation \cite{fleming2006controlled}:
\begin{equation}\label{HJBflem}
\frac{\partial}{\partial t}S(t,x) + H(t,x,D_{x}S(t,x), D_{xx}S(t,x))=0,
\end{equation}
where 
$H$ is some Hamiltonian. Thorough analysis of equations such as (\ref{HJBflem}) is presented in \cite{fleming2006controlled}. These HJB equations have been proven to have viscosity solutions, which form a particular class of weak solutions, and when we analyse our fractional HJB type equation it is reasonable to look for viscosity solutions too. This topic will also be discussed in the authors' further research.

\section{Fractional HJB for a controlled scaled CTRW with running rewards}

\smallskip

If we now take into account the reward for jumping from position $y$, namely $f(u,y)$, then the DP equation for the optimal payoff becomes the following: 
\begin{eqnarray}
S^{\tau}(t,y)=\sup_{u \in U} \left[[S^{\tau}_{0}(y)]\int_{t}^{\infty}\nu(dr/\tau^{1/\beta})+ \int_{0}^{t} \int_{\mathbb{R}^{d}} \left[S^{\tau}(t-r,y+\xi) + f(u,y,\xi) \right] \mu_{u}(d\xi/\tau^{1/\alpha}) \nu(dr/\tau^{1/\beta})\right].
\end{eqnarray} 

We re-write it in the form:
\begin{eqnarray}
S^{\tau}(t,y)= \sup_{u \in U}\left[S^{\tau}_{0}(y)\int_{t}^{\infty}\nu(dr/\tau^{1/\beta})+ \int_{0}^{t} \int_{\mathbb{R}^{d}} \left[S^{\tau}(t-r,y+\xi)\right] \mu_{u}(d\xi/\tau^{1/\alpha}) \nu(dr/\tau^{1/\beta})\right.+ \nonumber \\
\left. +\int_{0}^{t}\int_{\mathbb{R}^{d}} (f(u,y,\xi)) \mu_{u}(d\xi/\tau^{1/\alpha}) \nu(dr/\tau^{1/\beta})\right].
\end{eqnarray}

Now, we make an assumption on the function $f$ as follows: for $\lambda>0$ 
\begin{eqnarray}
f(u,y,\lambda \xi) = \lambda^{\alpha}f(u,y,\xi),
\end{eqnarray}
i.e. $f$ is homogeneous of order $\alpha$, where $\alpha \in (0,2]$ is the stability index we've met before. Denote
\begin{equation}
F(u,y)=\int_{\mathbb{R}^{d}}f(u,y,\xi)\mu_{u}(d\xi).
\end{equation}

We have for $\tau \rightarrow 0$:
\begin{eqnarray}
\frac{1}{\tau}\int_{0}^{t}\int_{\mathbb{R}^{d}}f(u,y,\xi \tau^{1/\alpha})\mu_{u}(d\xi)\nu(dr/\tau^{1/\beta})= 
\frac{1}{\tau}\int_{0}^{t\tau^{-1/\beta}}\int_{\mathbb{R}^{d}}\tau f(u,y,\xi)\mu_{u}(d\xi)\nu(dr)= \nonumber \\
F(u,y)\int_{0}^{t \tau^{-1/\beta}}\nu(dr) \rightarrow F(u,y).
\end{eqnarray}

Therefore, the full limiting equation with the included reward function for jumping, obtains the following form:

\begin{empheq}[box={\color{black}\fboxsep=10pt\fbox}]{align}
\color{red}
0= \frac{t^{-\beta}}{\Gamma(1-\beta)}\hat{S}_{0}(y) + \sup_{u \in U}\left[L_{\alpha}^{u}\hat{S}(t,y) +F(u,y) \right] +A^{*}_{\beta}\hat{S}(t,y).
\end{empheq}

Note that $A^{*}_{\beta}\hat{S}(t,y)$ is in independent of control, so it can be taken out of the supremum. If we also include waiting rewards into the system, the DP equation for $S^{\tau}(t,y)$ then becomes
\begin{eqnarray}
S^{\tau}(t,y)= \sup_{u \in U}\left[S^{\tau}_{0}(y)\int_{t}^{\infty}\nu(dr/\tau^{1/\beta})+ \int_{0}^{t} \int_{\mathbb{R}^{d}} \left[S^{\tau}(t-r,y+\xi)\right] \mu_{u}(d\xi/\tau^{1/\alpha}) \nu(dr/\tau^{1/\beta})\right. \nonumber \\
\left. +\int_{0}^{t}\int_{\mathbb{R}^{d}} (f(u,y,\xi)) \mu_{u}(d\xi/\tau^{1/\alpha}) \nu(dr/\tau^{1/\beta})\right.  \nonumber \\
+\left.  t\int_{t}^{\infty}g_{u,t}(y)\nu(dr/\tau^{\beta})+ \int_{0}^{t}\int_{\mathbb{R}^{d}} g_{u,t}(y)r\mu_{u}(d\xi/\tau^{1/\alpha})\nu(dr/\tau^{1/\beta}) \right]
\end{eqnarray}

and for $\tau \rightarrow 0$ we have
\begin{eqnarray}
\frac{g_{u,t}(y)}{\tau}\int_{0}^{t}\nu(dr/ \tau^{1/\beta})=
\frac{g_{u,t}(y)}{\tau}\int_{0}^{t \tau^{-1/\beta}}\frac{dr}{r^{1+\beta}}=
\frac{g_{u,t}(y)}{\tau}(t \tau^{-1/\beta})^{-\beta}=g_{u,t}(y)t^{-\beta},
\end{eqnarray}

and
\begin{eqnarray}
\frac{1}{\tau}\int_{0}^{t}\int_{\mathbb{R}^{d}} g_{u,t}(y)r\mu_{u}(d\xi/\tau^{1/\alpha})\nu(dr/\tau^{1/\beta}) =
\frac{1}{\tau}g_{u,t}(y)\int_{0}^{t \tau^{-1/\beta}} \frac{\tau^{1/\beta} r dr}{r^{1+\beta}} = \nonumber \\
g_{u,t}(y)\frac{\tau^{1/\beta}}{\tau}(t \tau^{-1/\beta})^{1-\beta}=
g_{u,t}(y)t^{1-\beta}.
\end{eqnarray}

This gives us:

\begin{empheq}[box={\color{black}\fboxsep=10pt\fbox}]{align}
\color{red}
0= \frac{t^{-\beta}}{\Gamma(1-\beta)}\hat{S}_{0}(y) + \sup_{u \in U}\left[L_{\alpha}^{u}\hat{S}(t,y) + 2g_{u,t}(y)t^{1-\beta} +F(u,y) \right] +A^{*}_{\beta}\hat{S}(t,y).
\end{empheq}

\section{Fractional HJB for a scaled CTRW with inner motion}

We can also add particle walking during waiting, which is an arbitrary Feller process $W^{u}_{t}(y)$ and which we will refer to as inner motion. In \cite{gihman1975theory} such a process, i.e. with inner motion, is referred to as {\it a process with a semi-Markov chance interference}. Our DP equation obtains the following form:

\begin{equation}
S^{\tau}(t,y)= \sup_{u \in U} \left[\left(\mathbb{E}[S^{\tau}(W^{u}_{t}(y)\right)P(\tau > t)]+ \int_{0}^{t}\int_{\mathbb{R}^{d}}\mathbb{E}\left(S^{\tau}(t-r,W_{t}^{u}(y)+\xi \right)| \gamma = r)\mu_{u}(d\xi/\tau^{1/\alpha})\nu(dr/\tau^{1/\beta})\right],
\end{equation}

where $W_{t}^{y}$ is a limiting Markov process and a generator $L$ for the limit of the scaled controlled CTRW, as $\tau \rightarrow 0$. We scale the waiting times for this inner process by $\tau^{1/\beta}$ too, and also
\begin{equation}
\mathbb{E}S(W_{r}(y)+\xi, t-r)=\int_{\mathbb{R}^{d}}S(z+\xi,t-r)P_{r}(y,dz),
\end{equation}
and
\begin{equation}
\mathbb{E}[S^{\tau}_{0}(W^{u}_{t}(y))]=\int_{\mathbb{R}^{d}}S^{\tau}_{0}(z)P_{t}(W^{u}_{t}(y),dz),
\end{equation}
and in the limit as $\tau \rightarrow 0$ the inner motion has a generator $B$. So we have the following: 

\begin{eqnarray}
S^{\tau}(t,y)=\sup_{u \in U}\left[\int_{t}^{\infty}\int_{\mathbb{R}^{d}}S^{\tau}(z)P_{t}(y,dz)\nu(dr/\tau^{1/\beta}) \right.  \nonumber \\ 
+\left. \int_{0}^{t}\int_{\mathbb{R}^{d}}\int_{\mathbb{R}^{d}}S^{\tau}(z+\xi,t-r)P_{r}(y,dz)\mu_{u}(d\xi/\tau^{1/\alpha})\nu(dr/\tau^{1/\beta})\right]\nonumber \\
=\sup_{u \in U}\left[\int_{t\tau^{-1/\beta}}^{\infty}\int_{\mathbb{R}^{d}}S^{\tau}(z)P_{t \tau^{1/\beta}}(y,dz)\nu(dr)  \right. \nonumber \\ 
+\left. \int_{0}^{t\tau^{-1/\beta}}\int_{\mathbb{R}^{d}}\int_{\mathbb{R}^{d}}S^{\tau}(z+\xi \tau^{1/\alpha},t-r\tau^{1/\beta})P_{r \tau^{1/\beta}}(y,dz)\mu_{u}(d\xi)\nu(dr)\right]\nonumber \\
=\sup_{u \in U}\left[
\int_{t\tau^{-1/\beta}}^{\infty}\int_{\mathbb{R}^{d}}S^{\tau}(z)P_{t \tau^{1/\beta}}(y,dz)\nu(dr) \right. \nonumber \\
 +\int_{0}^{t\tau^{-1/\beta}}\tau \int_{\mathbb{R}^{d}}\frac{\int_{\mathbb{R}^{d}}
S^{\tau}(z+\xi \tau^{1/\alpha},t-r\tau^{1/\beta})P_{r \tau^{1/\beta}}(y,dz)-
 S^{\tau}(y,t-r\tau^{1/\beta})}{\tau}\mu_{u}(d\xi)\nu(dr)   \nonumber \\
\left. + \int_{0}^{t\tau^{-1/\beta}} \left(\tau \frac{S^{\tau}(y,t-r\tau^{1/\beta})- S^{\tau}(t,y)}{\tau}+S^{\tau}(t,y)\right)\nu(dr)\right].
\end{eqnarray}

Hence,

\begin{eqnarray}
\frac{1}{\tau} \left[ \int_{t\tau^{-1/\beta}}^{\infty}S^{\tau}(t,y)\nu(dr)- \int_{t\tau^{-1/\beta}}^{\infty}\int_{\mathbb{R}^{d}}S^{\tau}(z)P_{t\tau^{1/\beta}}(y,dz)\nu(dr)\right] \nonumber \\
=\sup_{u \in U}\left[\int_{0}^{t\tau^{-1/\beta}} \frac{S^{\tau}(t-r\tau^{1/\beta},y)- S^{\tau}(t,y)}{\tau} \nu(dr) \right. \nonumber \\ 
+\int_{0}^{t\tau^{-1/\beta}} \int_{\mathbb{R}^{d}}\int_{\mathbb{R}^{d}}\frac{
S^{\tau}(t-r\tau^{1/\beta},z+\xi \tau^{1/\alpha})- S^{\tau}(t-r\tau^{1/\beta},z)}{\tau}P_{r\tau^{1/\beta}}(y,dz)\mu_{u}(d\xi)\nu(dr)\nonumber \\
+\left. \int_{0}^{t\tau^{-1/\beta}} \int_{\mathbb{R}^{d}}\frac{S^{\tau}(t-r\tau^{1/\beta},z)-S^{\tau}(t-r\tau^{1/\beta},y)}{\tau}P_{r\tau^{1/\beta}}(y,dz)\nu(dr)\right].
\end{eqnarray}

The last summand can be represented in the form:
\begin{eqnarray}
\int_{0}^{t\tau^{-1/\beta}} \int_{\mathbb{R}^{d}}\frac{S^{\tau}(t-r\tau^{1/\beta},z)-S^{\tau}(t-r\tau^{1/\beta},y)}{\tau}P_{r\tau^{1/\beta}}(y,dz)\nu(dr) \nonumber \\
=\int_{0}^{\infty} \int_{\mathbb{R}^{d}}\frac{S^{\tau}(t-r\tau^{1/\beta},z)-S^{\tau}(t-r\tau^{1/\beta},y)}{\tau}P_{r\tau^{1/\beta}}(y,dz)\nonumber \\
-\int_{t\tau^{-1/\beta}}^{\infty} \int_{\mathbb{R}^{d}}\frac{S^{\tau}(t-r\tau^{1/\beta},z)-S^{\tau}(t-r\tau^{1/\beta},y)}{\tau}P_{r\tau^{1/\beta}}(y,dz) \nu(dr).
\end{eqnarray}

As $\tau \rightarrow 0$, due to the boundary condition the two terms in the second integral are $0$'s and the first integral satisfies:
\begin{eqnarray}
\int_{0}^{\infty} \int_{\mathbb{R}^{d}}\frac{S^{\tau}(t-r\tau^{1/\beta},z)-S^{\tau}(t-r\tau^{1/\beta},y)}{\tau}P_{r\tau^{1/\beta}}(y,dz) \nu(dr) \rightarrow B^{u}\hat{S}(t,y)
\end{eqnarray} 

As for the term 
\begin{equation}
\frac{1}{\tau}\int_{t\tau^{-1/\beta}}^{\infty}\int_{\mathbb{R}^{d}}S^{\tau}(z)P_{t \tau^{1/\beta}}(y,dz)\nu(dr), 
\end{equation}
we add and subtract 
\begin{equation}
\frac{1}{\tau}\int_{t \tau^{-1/\beta}}^{\infty}S^{\tau}_{0}(y)\nu(dr)
\end{equation}
and this gives us:
\begin{eqnarray}
\frac{1}{\tau}\int_{t\tau^{-1/\beta}}^{\infty}\int_{\mathbb{R}^{d}}S^{\tau}(z)P_{t \tau^{1/\beta}}(y,dz)\nu(dr)\nonumber \\
=\int_{t\tau^{-1/\beta}}^{\infty}\frac{\left[\int_{\mathbb{R}^{d}}S^{\tau}(z)P_{t \tau^{1/\beta}}(y,dz) - S^{\tau}_{0}(y)\right]}{\tau} \nu(dr) + \frac{1}{\tau}\int_{t\tau^{-1/\beta}}^{\infty}S^{\tau}_{0}(y)\nu(dr).
\end{eqnarray}

Now,
\begin{eqnarray}
\int_{t \tau^{-1/\beta}}^{\infty}t \tau^{1/\beta - 1}\int_{\mathbb{R}^{d}}\frac{S^{\tau}(z)- S^{\tau}_{0}(y)}{t\tau^{1/\beta}}P_{t \tau^{1/\beta}}(y,dz)\nu(dr) \rightarrow \int_{t\tau^{-1/\beta}}^{\infty}t\tau^{1/\beta - 1}B^{u}S^{\tau}_{0}(y)\nu(dr)
\rightarrow 0,
\end{eqnarray}

and
\begin{equation}
\frac{1}{\tau}\int_{t\tau^{-1/\beta}}^{\infty}S^{\tau}_{0}(y)\nu(dr) \sim \frac{t^{-\beta}}{\Gamma(1-\beta)}\hat{S}_{0}(y),
\end{equation}

Hence the equation for $\hat{S}(t,y)$ with an inner motion during waiting, as $\tau \rightarrow 0$, is 

\begin{empheq}[box={\color{black}\fboxsep=10pt\fbox}]{align}
\color{red}
0=\frac{1}{\Gamma(1-\beta)}t^{-\beta}\hat{S}_{0}(y)+ \sup_{u \in U}\left[L^{u}_{\alpha}\hat{S}(t,z)+B^{u}\hat{S}(t,y)\right]+A^{*}_{\beta}\hat{S}(t,y).
\end{empheq}

If we also add the rewards for jumping and rewards for waiting, then the extra terms with $g_{u,t}(y)$ and $F(u,y)$ appear, yielding:

\begin{empheq}[box={\color{black}\fboxsep=30pt\fbox}]{align}
\color{red}
0=\frac{1}{\Gamma(1-\beta)}t^{-\beta}\hat{S}_{0}(y)+ \sup_{u \in U}\left[L^{u}_{\alpha}\hat{S}(t,z)+B^{u}\hat{S}(t,y) + 2g_{u,t}(y)t^{1-\beta} +F(u,y) \right]+A^{*}_{\beta}\hat{S}(t,y).
\end{empheq}

\section{The HJB written backwards in time}

When working with time non-homogeneous processes it is more convenient to change the notation.
Now $S(t,y)$ will denote the optimal pay-off at position $y$ and time $t$, where $t$ is the time elapsed since the beginning of the process. Now we will have $S(T,y)$ known, where $T$ is the terminal time of the process. Now $S(t,y)$ is defined as

\begin{equation}
S(t,y)=\sup_{\tilde{u} \in \tilde{U}}\mathbb{E}[S(T,y+Y(T-t,\tilde{u}))].
\end{equation}

We can write the DP equation backwards, which means time is moving forwards and $t$ denotes time elapsed since the beginning of the process. We shall see the differences this brings to the DP equation as $\tau \rightarrow 0$ is that firstly, instead of the dual operator $A^{*}$ we will have the generator $A$ and secondly, that the boundary term will depend on $T$ too.

\begin{eqnarray}
S^{\tau}(t,y)=\sup_{u \in U}\left[S^{\tau}(T,y)\int_{(T-t)\tau^{-1/\beta}}^{\infty}\nu(dr) \right.  \nonumber \\ 
+\left. \tau \int_{0}^{(T-t)\tau^{-1/\beta}}\int_{\mathbb{R}^{d}}\frac{S^{\tau}(t+r\tau^{1/\beta},y+\xi \tau^{1/\alpha})-S^{\tau}(t+r \tau^{1/\beta},y)}{\tau} \mu_{u}(d\xi)\nu(dr) \right.  \nonumber \\
+\left. \tau \int_{0}^{(T-t)\tau^{-1/\beta}}\int_{\mathbb{R}^{d}}\frac{S^{\tau}(t+r\tau^{1/\beta},y)-S^{\tau}(t,y)}{\tau}\mu_{u}(d\xi)\nu(dr) \right.  \nonumber \\
+\left.  \int_{0}^{(T-t)\tau^{-1/\beta}}\int_{\mathbb{R}^{d}}S^{\tau}(t,y)\mu_{u}(d\xi)\nu(dr) \right].
\end{eqnarray}

We divide the whole equation by $\tau$. As $\tau \rightarrow 0$ the first term on the RHS gives us the limiting boundary term 
\begin{equation}
\frac{(T-t)^{-\beta}}{\Gamma(1-\beta)}\hat{S}(T,y).
\end{equation}
Note that it is not dependent on control and hence we can take it out of the supremum.

The second term on the RHS gives us, as $\tau \rightarrow 0$, 
\begin{eqnarray}
\int_{0}^{(T-t)\tau^{-1/\beta}}\int_{\mathbb{R}^{d}}\frac{[S^{\tau}(t+r\tau^{1/\beta},y+\xi \tau^{1/\alpha})-S^{\tau}(t+r\tau^{1/\beta},y)]}{\tau} \mu(d\xi)\nu(dr) \rightarrow L_{\alpha}\hat{S}(t,y).
\end{eqnarray}

The third term we write in the following form:
\begin{eqnarray}
\int_{0}^{(T-t)\tau^{-1/\beta}}\int_{\mathbb{R}^{d}}\frac{[S^{\tau}(t+r\tau^{1/\beta},y)-S^{\tau}(t,y)]}{\tau}\mu_{u}(d\xi)\nu(dr)  \nonumber \\
=\int_{0}^{\infty}\int_{\mathbb{R}^{d}}\frac{[S^{\tau}(t+r\tau^{1/\beta},y)-S^{\tau}(t,y)]}{\tau}\mu_{u}(d\xi)\nu(dr)  \nonumber \\
- \int_{(T-t)\tau^{-1/\beta}}^{\infty}\int_{\mathbb{R}^{d}}\frac{[S^{\tau}(t+r\tau^{1/\beta},y)-S^{\tau}(t,y)]}{\tau}\mu_{u}(d\xi)\nu(dr) 
\end{eqnarray}

For the first integral above:
\begin{equation}
\int_{0}^{\infty}\int_{\mathbb{R}^{d}}\frac{[S^{\tau}(t+r\tau^{1/\beta},y)-S^{\tau}(t,y)]}{\tau}\mu_{u}(d\xi)\nu(dr)  \rightarrow A_{\beta}\hat{S}(t,y).
\end{equation}

As for the term

\begin{equation}
\int_{(T-t)\tau^{-1/\beta}}^{\infty}\int_{\mathbb{R}^{d}}\frac{S^{\tau}(t+r\tau^{1/\beta},y)}{\tau}\mu_{u}(d\xi)\nu(dr)=0 \mbox{ due to the boundary condition $S^{\tau}(t,y)=0$ for $t>T$,} 
\end{equation}

and
\begin{eqnarray}
\frac{1}{\tau}\int_{0}^{(T-t)\tau^{-1/\beta}}\int_{\mathbb{R}^{d}}S^{\tau}(t,y)\mu_{u}(d\xi)\nu(dr)\nonumber \\
=\frac{1}{\tau}\int_{0}^{\infty}\int_{\mathbb{R}^{d}}S^{\tau}(t,y)\mu_{u}(d\xi)\nu(dr)-\frac{1}{\tau}\int_{(T-t)\tau^{-1/\beta}}^{\infty}\int_{\mathbb{R}^{d}}S^{\tau}(t,y)\mu_{u}(d\xi)\nu(dr).
\end{eqnarray}

The first integral in the expression above for the last term of the RHS cancels out with $S^{\tau}(t,y)$ on the LHS, $\frac{1}{\tau}S^{\tau}(t,y)$ and the second term in the above representation cancels out with the term
\begin{equation}
-\int_{(T-t)\tau^{-1/\beta}}^{\infty}\int_{\mathbb{R}^{d}}\frac{-S^{\tau}(t,y)}{\tau}\mu_{u}(d\xi)\nu(dr). 
\end{equation}

Therefore we arrive at the forward equation as $\tau \rightarrow 0$:
\begin{empheq}[box={\color{black}\fboxsep=10pt\fbox}]{align}
\color{red}
0=\frac{(T-t)^{-\beta}\hat{S}(T,y)}{\Gamma(1-\beta)}+\sup_{u \in U}[L_{\alpha}^{u}\hat{S}(t,y)]+A_{\beta}\hat{S}(t,y).
\end{empheq} 

If jumping and waiting rewards are present, then we obtain:

\begin{empheq}[box={\color{black}\fboxsep=10pt\fbox}]{align}
\color{red}
0=\frac{(T-t)^{-\beta}\hat{S}(T,y)}{\Gamma(1-\beta)}+\sup_{u \in U}\left[L^{u}_{\alpha}\hat{S}(t,y) + 2g_{u,t}(y)t^{1-\beta} + F(u,y) \right]+A_{\beta}\hat{S}(t,y).
\end{empheq} 

If the inner motion is present too:

\begin{empheq}[box={\color{black}\fboxsep=10pt\fbox}]{align}
\color{red}
0=\frac{(T-t)^{-\beta}\hat{S}(T,y)}{\Gamma(1-\beta)}+\sup_{u \in U}\left[L^{u}_{\alpha}\hat{S}(t,y)+B^{u}\hat{S}(t,y) + 2g_{u,t}(y)t^{1-\beta} + F(u,y) \right]+A_{\beta}\hat{S}(t,y).
\end{empheq}

Here $\leftidx{_{\infty-}}{D^{\beta}}{_x}=A_{\beta}$. The proof is just below and consists of changing variables $x-y=z$, differentiating under the integral and using that $-\beta\Gamma{(-\beta)}=\Gamma{(1-\beta)}$ for $\beta \in (0,1)$, and changing variables again: $-z=y$:

\begin{eqnarray}
\leftidx{_{\infty-}}{D^{\beta}}{_x}f(x)=\frac{1}{\Gamma{(1-\beta)}}\frac{d}{dx}\int_{x}^{\infty}f(y)(y-x)^{-\beta}dy\nonumber \\
=\frac{1}{\Gamma{(1-\beta)}}\frac{d}{dx}\int_{-\infty}^{0}f(x-z)|z|^{-\beta}dz\nonumber \\
=\frac{1}{\Gamma{(1-\beta)}}\int_{-\infty}^{0}-\frac{\partial}{\partial z}(f(x-z)-f(x))|z|^{-\beta}dz\nonumber \\
=-\frac{1}{\Gamma{(1-\beta)}}|z|^{-\beta}(f(x-z)-f(x))|_{-\infty}^{0}\nonumber \\
+ \frac{1}{\Gamma{(1-\beta)}}\int_{-\infty}^{0}\frac{1}{-\beta}(f(x-z)-f(x))|z|^{-1-\beta}dz \nonumber \\
=-\frac{1}{\Gamma{(-\beta)}}\int_{-\infty}^{0}(f(x-z)-f(x))|z|^{-\beta-1}dz \nonumber \\
=\frac{1}{\Gamma{(-\beta)}}\int_{\infty}^{0}(f(x+y)-f(x))|y|^{-\beta-1}dy \nonumber \\
=-\frac{1}{\Gamma{(-\beta)}}\int_{0}^{\infty}(f(x+y)-f(x))|y|^{-\beta-1}dy \nonumber \\
=A_{\beta}f(x).
\end{eqnarray}

The boundary term appearing after integration by parts is $0$ because $|z|^{-\beta} \rightarrow 0$ as $|z| \rightarrow \infty$.

\section{Position and time dependent HJB}

\subsection{Position dependent HJB}

It is known that for an arbitrary Feller process $\hat{Y}_{t}$ with a generator $L^{u}$ it is always possible to construct a family of measures $\mu_{u,\tau, n},$ such that
\begin{equation}
Y_{n}^{\tau}=Y_{n-1}^{\tau}+\xi_{n}^{\tau}(Y_{n-1}^{\tau}),
\end{equation}
where $\xi_{n}, n \in \mathbb{N}$ are r.v.'s with distributions given by $\mu_{u,\tau,n}$, so that
\begin{equation}
\int_{\mathbb{R}^{d}}\frac{f(t,y+\xi)-f(t,y)}{\tau}\mu_{u,\tau}(y,d\xi) \rightarrow L^{u}f(t,y).
\end{equation}

Consequently for $n=[t/\tau]$, by theorem $19.25$, attributed to Trotter, Sova, Kurtz and Mackevicius, and theorem $19.28$ in \cite{kallenberg_foundations},

\begin{equation}
Y_{[t/\tau]}^{\tau} \rightarrow \hat{Y}_{t}.
\end{equation}

The DP equation is
\begin{eqnarray}
S^{\tau}(t,y)=\sup_{u \in U}\left(S^{\tau}_{0}(y)\int_{t}^{\infty}\nu(dr/\tau^{1/\beta})+
\int_{0}^{t} \int_{\mathbb{R}^{d}} S^{\tau}(t-r,y+\xi)\mu_{u,y,\tau}(d\xi) \nu(dr/\tau^{1/\beta})\right).
\end{eqnarray}
In the limit $\tau \rightarrow 0$ this turns into

\begin{eqnarray}
0=\frac{1}{\Gamma(1-\beta)} t^{-\beta}\hat{S}_{0}(y)+\sup_{u \in U}[L^{u}\hat{S}(t,y)] + A^{*}_{\beta}\hat{S}(t,y).
\end{eqnarray}

If the inner motion and jumping and waiting rewards are also included, we have:
\begin{empheq}[box={\color{black}\fboxsep=10pt\fbox}]{align}
\color{red}
0=\frac{1}{\Gamma(1-\beta)} t^{-\beta}\hat{S}_{0}(y)+\sup_{u \in U}\left[L^{u}\hat{S}(t,y)+B^{u}\hat{S}(t,y) + 2g_{u,t}(y)t^{1-\beta} + F(u,y) \right] + A^{*}_{\beta}\hat{S}(t,y).
\end{empheq} 

Here the generator $L^{u}$ is an arbitrary generator of a Feller process $\hat{Y}_{t}$ and $A_{\beta}^{*}$ is the dual operator for a generator $A$ of a $\beta$-stable Feller process $\hat{X}_{t}$, and $B$ is the generator of the Feller process which is the inner motion.

\subsection{Time dependent HJB}

We could have stable position dependence together with general time dependence.
In this case

\begin{equation}
X_{n}^{\tau}=X_{n-1}^{\tau}+\gamma_{n}^{\tau}(X_{n-1}^{\tau}).
\end{equation}

We build families of $\gamma_{n}^{\tau}$, $n=[t/\tau] \in \mathbb{N}$, each with the law $\nu_{t,\tau}(X_{n-1}^{\tau},dr)$ so that we can have a discrete time Markov chain approximation of an arbitrary increasing Feller process $\hat{X}_{t}$ with a generator $A$:
\begin{equation}
\int_{0}^{\infty}\frac{f(t+r)-f(t)}{\tau}\nu_{\tau}(t,dr) \rightarrow Af(t).
\end{equation}
Consequently by theorems $19.25$ and $19.28$ in \cite{kallenberg_foundations}
\begin{equation}
X_{[t/\tau]}^{\tau} \rightarrow \hat{X}_{t}.
\end{equation}
For instance, in the case when $A$ is stable-like we have:

\begin{equation}
\int_{t\tau^{-1/\beta(t)}}^{\infty}\nu(dr) \sim \frac{\zeta(t)}{\Gamma(1-\beta(t)) (t\tau^{-1/\beta(t)})^{\beta(t)}}
\end{equation}
as $\tau \rightarrow 0$ and so $t \tau^{-1/\beta} \rightarrow \infty$. Note that now the measure $\nu$ depends on time $t$ elapsed since the beginning of the process. Also, we include waiting time rewards. Now boundary condition is $S(T,y)$ is known and we write the DP equation:
\begin{eqnarray}
S^{\tau}(t,y)=\sup_{u \in U}\left[\left[S^{\tau}(T,y)+g_{u,t}(y)(T-t)\right]\int_{T-t}^{\infty}\nu_{t,\tau}(dr) \right. \nonumber \\
\left. + \int_{0}^{T-t}\int_{\mathbb{R}^{d}}\left[S^{\tau}(t+r,y+\xi)+g_{u,t}(y)r + f(y+\xi,u)\right] \mu_{u, y}(d\xi/ \tau^{1/\alpha})\nu_{t,\tau}(dr)\right].
\end{eqnarray}

We assume that
\begin{equation}
\lim_{\tau \rightarrow 0}\frac{1}{\tau}\int_{0}^{t}r \nu_{\tau}(dr)=\kappa_{1}(t),
\end{equation}
and
\begin{equation}
\lim_{\tau \rightarrow 0}\frac{1}{\tau}\int_{t}^{\infty} \nu_{\tau}(dr)=\kappa_{2}(t),
\end{equation}
and that
\begin{equation}
\frac{1}{\tau}\int_{\mathbb{R}^{d}}f(u,y,\xi)\mu_{u}(d \xi / \tau^{1/\alpha})=F(u,y).
\end{equation}
 
This yields us the following limiting equation  for $\hat{S}(t,y)$ with an arbitrary generator of an increasing Feller process $\hat{X}_{t}$, namely $A$, and a generator of a stable Feller process $\hat{Y}_{t}$, namely $L_{\alpha}^{u}$, as $\tau \rightarrow 0$:

\begin{empheq}[box={\color{black}\fboxsep=10pt\fbox}]{align}
\color{red}
0=A\mathbf{1}_{(T-t>0)} \hat{S}(T,y)+A\hat{S}(t,y) \nonumber \\
\color{red}
+\sup_{u \in U}\left[ L_{\alpha}^{u}\hat{S}(t,y)+g_{u,t}(y)(T-t)A\mathbf{1}_{(T-t>0)}\right. \nonumber \\
\color{red}
\left.+ g_{u,t}(y)(Ac)(t)  + F(u,y)A\mathbf{1}_{(T-t>0)}\right],
\end{empheq}
where $c(t)=t\mathbf{1_{(T-t>0)}}$.

In case when inner motion also takes place, we have:
\begin{empheq}[box={\color{black}\fboxsep=10pt\fbox}]{align}
\color{red}
0=A\mathbf{1}_{(T-t>0)} \hat{S}(T,y)+A\hat{S}(t,y) \nonumber \\
\color{red}
+\sup_{u \in U}\left[ L_{\alpha}^{u}\hat{S}(t,y)+ B^{u}\hat{S}(t,y) + g_{u,t}(y)(T-t)A\mathbf{1}_{(T-t>0)}\right.  \nonumber \\
\color{red}
\left.+ g_{u,t}(y)(Ac)(t) + F(u,y)A\mathbf{1}_{(T-t>0)}\right],
\end{empheq}
where $c(t)=t\mathbf{1_{(T-t>0)}}$, for $t \ge 0$. Note that the key difference is that there are no  $\beta$'s involved at all in this case, and the generator $A$ is an arbitrary one describing an increasing Feller process $\hat{X}_{t}$.

\subsection{Both time and position dependent HJB}

We could have position dependence together with general time dependence:

\begin{equation}
X_{n}^{\tau}=X_{n-1}^{\tau}+\gamma_{n}^{\tau}(X_{n-1}^{\tau}),
\end{equation}

\begin{equation}
Y_{n}^{\tau}=Y_{n-1}^{\tau}+\xi_{n}^{\tau}(Y_{n-1}^{\tau}).
\end{equation}

We build families of $\gamma_{i}^{\tau}$ and $\xi_{i}^{\tau}$ as before. Also, we include waiting time rewards into the system. Again we use theorems $19.25$ and $19.28$ in \cite{kallenberg_foundations} which guarantee that such approximations of the Feller processes by discrete time Markov chains are possible. Our boundary condition is $\hat{S}(T,y)$ and we write the DP equation:
\begin{eqnarray}
S^{\tau}(t,y)=\sup_{u \in U}\left[\left[S^{\tau}(T,y)+g_{u,t}(y)(T-t)\right]\int_{T-t}^{\infty}\nu_{t,\tau}(dr) \right. \nonumber \\
\left. + \int_{0}^{T-t}\int_{\mathbb{R}^{d}}\left[S^{\tau}(t+r,y+\xi)+g_{u,t}(y)r + f(y+\xi,u)\right] \mu_{u,\tau, y}(d\xi)\nu_{t,\tau}(dr)\right].
\end{eqnarray}

For the DP equation not to blow up, we assume that
\begin{equation}
\lim_{\tau \rightarrow 0}\frac{1}{\tau}\int_{0}^{t}r \nu_{\tau}(dr)=\kappa_{1}(t),
\end{equation}
and
\begin{equation}
\lim_{\tau \rightarrow 0}\frac{1}{\tau}\int_{t}^{\infty} \nu_{\tau}(dr)=\kappa_{2}(t),
\end{equation}
and that
\begin{equation}
\frac{1}{\tau}\int_{\mathbb{R}^{d}}f(u,y,\xi)\mu_{u,\tau}(d \xi)=F(u,y).
\end{equation}
 
This yields us the following limiting equation with arbitrary generators of Feller processes $\hat{X}_{t}$ and $\hat{Y}_{t},$ namely $A$ and $L$, as $\tau \rightarrow 0$:

\begin{empheq}[box={\color{black}\fboxsep=10pt\fbox}]{align}
\color{red}
0=A\mathbf{1}_{(T-t>0)} \hat{S}(T,y)+A\hat{S}(t,y) \nonumber \\
\color{red}
+\sup_{u \in U}\left[ L\hat{S}(t,y)+g_{u,t}(y)(T-t)A\mathbf{1}_{(T-t>0)}\right. \nonumber \\
\color{red}
\left.+ g_{u,t}(y)(Ac)(t)  + F(u,y)A\mathbf{1}_{(T-t>0)}\right],
\end{empheq}
where $c(t)=t\mathbf{1_{(T-t>0)}}$.

In case when inner motion also takes place, we have:
\begin{empheq}[box={\color{black}\fboxsep=10pt\fbox}]{align}
\color{red}
0=A\mathbf{1}_{(T-t>0)} \hat{S}(T,y)+A\hat{S}(t,y) \nonumber \\
\color{red}
+\sup_{u \in U}\left[ L\hat{S}(t,y)+ B^{u}\hat{S}(t,y) + g_{u,t}(y)(T-t)A\mathbf{1}_{(T-t>0)}\right.  \nonumber \\
\color{red}
\left.+ g_{u,t}(y)(Ac)(t) + F(u,y)A\mathbf{1}_{(T-t>0)}\right],
\end{empheq}
where $c(t)=t\mathbf{1_{(T-t>0)}}$, for $t \ge 0$. Note that the key difference is that now there are no $\alpha$'s and $\beta$'s involved at all in the equation for $\hat{S}(t,y)$, i.e. the generators of the increasing Feller process $\hat{X}_{t}$ and a Feller process $\hat{Y}_{t}$ are arbitrary.

\section{Further research}
In our further research is proving uniqueness and existence of classical solutions to fractional HJB type equations with Caputo fractional derivatives, and also proving convergence of $S^{\tau}(t,y)$ to $\hat{S}(t,y)$ as $\tau \rightarrow 0$ which corresponds precisely to the limiting process. Further work includes direct applications of this theory to finances \cite{kolokoltsov_gametheory}, \cite{kolokoltsov_gamefinances}.

\end{document}